\documentstyle[12pt,leqno,fleqn]{article}
\title{
\author{\bf{Lech Pasicki}}
\bf{A Stokes theorem for everyone}
}
\date{}
\newtheorem{theorem}{\indent Theorem}
\newtheorem{lemma}[theorem]{\indent Lemma}
\newtheorem{proposition}[theorem]{\indent Proposition}
\newtheorem{definition}[theorem]{\indent Definition}
\newtheorem{corollary}[theorem]{\indent Corollary}
\newtheorem{remark}[theorem]{\indent Remark}

\newcommand{\Int}{\mbox{Int\,}}
\newcommand{\Fr}{\mbox{Fr\,}}
\newcommand{\mx}{\mbox{max\,}}
\newcommand{\ckj}{\mbox{$c_{\stackrel{-}{j}}$}}
\newcommand{\tkj}{\mbox{$t_{\stackrel{-}{j}}$}}
\newcommand{\xkj}{\mbox{$x_{\stackrel{-}{j}}$}}
\newcommand{\Tkj}{\mbox{$T_{\stackrel{-}{j}}$}}

\newcommand{\Pkj}{\mbox{$P_{\stackrel{-}{j}}$}}
\newcommand{\deteIj}{\mbox{$(det \frac{De_{I}}{Dt_{\stackrel{-}{j}}})_{\mid j}$}}
\newcommand{\PhA}{\mbox{$\Phi\!: A \rightarrow R^{n}$}}
\newcommand{\PsA}{\mbox{$\Psi\!: A \rightarrow R^{n}$}}

\begin{document}
\maketitle
\vspace{1 in}

\begin{abstract}
Many versions of the Stokes theorem are known. More advanced of them require complicated
mathematical machinery to be formulated which discourages the users. Our theorem is sufficiently simple to
suit the handbooks and yet it is pretty general, as we assume the differential form to be continuous
on a compact set $  \Phi(A) $ and $ C^{1} $ ``inside'' while $  \Phi(A) $ is built of ``bricks'' and its inner
part is a $ C^{1} $ manifold. There is no problem of orientability and the integrals under consideration are
convergent. The proof is based on integration by parts and inner approximation.
\end{abstract}
First let us formulate our main result.
\begin{theorem}
\label{Th1}
Let $ A \subset R^{k} $, $(k>1) $ be a regular set (Def. \ref{De3}) consisting of normal
sets $ A_{1},\ldots,A_{m} $ and let $  \PhA  $ be a one to one 
mapping of class $ C^{1}$.  Assume that  $ \stackrel{k-1}{\omega}=\sum b_{I} \bigwedge dy_{I} $
is a continuous form on $ S=\Phi(A)$ and of class $ C^{1} $ on
$  \bigcup \{\Phi(\Int A_{j}): j=1,\ldots,m\}$. Then the integrals in the following formula are convergent
\begin{eqnarray*}
\qquad \int_{S}d\stackrel{k-1}{\omega}=
\int_{\partial S}\stackrel{k-1}{\omega} \mbox{,}
\end{eqnarray*}
and the equality holds on condition that the orientations of $ \partial S = \Phi(\Fr \Int A) $ and
of $S$ coincide.
\end{theorem}
\par Let us recall that the orientations under consideration coincide if
the orientation of the base of the hyperplane tangent to $S$ is the same as
that of the normal versor to $ \partial S $ being tangent and external
to $S$, supplemented with the base of vectors tangent to $ \partial S$, in points where
the normal versor under consideration exists.
\par It should be stressed that the integrals in Theorem \ref{Th1} do not
depend on parametrizations and that the assumptions guarantee the
orientability of $ S,\partial S$.
\begin{definition}
\label{De2}
  A set $ A \subset R^{k} $ is normal if its points, likely after 
renumbering of variables and mappings, satisfy
\begin{eqnarray}
\label{co1}
\left\{\begin{array}{l}
f_{1} \leq x_{1} \leq g_{1},\\
f_{2}(x_{1}) \leq x_{2} \leq g_{2}(x_{1}),\\
\vdots\\
f_{k}(x_{1},\ldots,x_{k-1}) \leq x_{k} \leq g_{k}(x_{1},\ldots,x_{k-1}),
\end{array}\right.
\end{eqnarray}
where $ f_{1},g_{1} $ are constants, all $ f_{i},g_{i} $ are continuous and
all $ f_{i\mid j},g_{i\mid j} $ are continuous and bounded inside the domains of 
$ f_{i},g_{i}$.
\end{definition}
\begin{definition}
\label{De3}
  A set $ A \subset R^{k} $ is regular if it is a finite sum of  normal sets
with disjoint interiors.
\end{definition}
\par For the sake of completeness we present the following well known basic
definitions and corollaries.
\begin{definition}
\label{De4}
   A k-form of the class $ C^{r} $ in $ R^{n} $ is an expression
\begin{eqnarray}
\label{co2}
  \stackrel{k}{\omega}
    =\sum_{i_{1},\ldots,i_{k}=1}^{n}\stackrel{k}{\omega}_{i_{1},\ldots,i_{k}}
    =\sum_{i_{1},\ldots,i_{k}=1}^{n}b_{i_{1},\ldots,i_{k}}(y)
    dy_{i_{1}}\wedge\cdots\wedge dy_{i_{k}},
\end{eqnarray}
where $  b_{i_{1},\ldots,i_{k}}\!: V \rightarrow R  $ are maps of class $ C^{r} $
on a set $  V \subset R^{n}$.
\end{definition}
  \par If $  A \subset R^{k} $ is a measurable set,$ \PhA  $
is a mapping and $  S=\Phi(A) \subset V$, we adopt
\begin{eqnarray}
\label{co3}
  \int_{S}\stackrel{k}{\omega}\mid_{\Phi}
  =\int_{A}\sum_{i_{1},\ldots,i_{k}=1}^{n}
  b_{i_{1},\ldots,i_{k}}(\Phi(x))
  det\frac{D(\Phi_{i_{1}},\ldots,\Phi_{i_{k}})}{D(x_{1},\ldots,x_{k})}dx
\end{eqnarray}
on condition that the integral on the right side has sense.
\par From the properties of determinants we obtain
\begin{corollary}
\label{Cor5}
  In formula (\ref{co2}) the components with $  i_{j}=i_{m} $ can be
disregarded, as the respective determinant equals zero.
\end{corollary}
\begin{corollary}
\label{Cor6}
  For $  k > n  $ we have $  \stackrel{k}{\omega}=0$.
\end{corollary}
\begin{corollary}
\label{Cor7}
  The transposition of $i_{j}$ and $i_{m}$ in
$  dy_{i_{1}}\wedge\cdots\wedge dy_{i_{k}}  $ changes the sign of
this expression.
\end{corollary}
\par We adopt in addition that 0-forms are mappings.
\begin{definition}
\label{De8}
   The exterior differential of a $  \stackrel{k-1}{\omega}$,
($k \geq 1$) is the following k-form
\begin{eqnarray}
\label{co4}
  d\stackrel{k-1}{\omega}
    =\sum_{i_{1},\ldots,i_{k-1}=1}^{n}db_{i_{1},\ldots,i_{k-1}}(y)\wedge
    dy_{i_{1}}\wedge\cdots\wedge dy_{i_{k-1}}=\\ \nonumber
    \sum_{i_{1},\ldots,i_{k-1}=1}^{n} \sum_{i=1}^{n}
     b_{i_{1},\ldots,i_{k-1}\mid i}(y)
    dy_{i}\wedge dy_{i_{1}}\wedge\cdots\wedge dy_{i_{k-1}}.
\end{eqnarray}
\end{definition}
\par For the mappings of class $ C^{2} $ the mixed second order partial derivatives are
equal and hence (Corollary \ref{Cor7}) we obtain
\begin{corollary}
\label{Cor9}
  IF $  \stackrel{k-1}{\omega}$, ($k \geq 1$) is a form of class $ C^{2} $
then we have $  d(d \stackrel{k-1}{\omega})=0$.
\end{corollary}
\par All functions in this paper are real valued. As regards $  t=(t_{1},\ldots,t_{k}) \in R^{k} $
we adopt $ \tkj=(t_{1},\ldots,t_{j-1},t_{j+1},\ldots,t_{k}) $ and 
$ (s,\tkj)=(t_{1},\ldots,t_{j-1},s,t_{j+1},\ldots,t_{k}) $ - in particular
$ (t_{j},t_{\stackrel{-}{j}})=t$.  For a segment 
$ T=<a_{1},b_{1}>\times \ldots \times <a_{k},b_{k}>\subset R^{k} $ we adopt 
$ T_{j}=<a_{j},b_{j}> $, $ [f]_{T_{j}}=f(b_{j})-f(a_{j}) $ and 
$ \Tkj=\{t_{\stackrel{-}{j}}\!: t \in T\}$.   The notation $ | | $ means
the modulus or the outer measure (of the respective dimension) which, as it is known, for 
compacta coincides with the Lebesgue measure.
\begin{proposition}[integration by parts]
\label{Pro10} 
  Let $ f,g $ be continuous on $ T=<a,b> $ and let
$ f',g' $ exist and be continuous on $ \Int T  $ with $ g' $ bounded .
Then we have
\begin{eqnarray}
\label{co5}
 \int_{T}f'g=[fg]_{T} - \int_{T}fg'.
\end{eqnarray}
\end{proposition}
\par Proof.  Let us adopt $ Q=<p,q> \subset (a,b)$.  Then the classical
theorem on integration by parts implies  $ \int_{Q}f'g=[fg]_{Q} - \int_{Q}fg'$.
Now by the continuity of  $ fg $ and the fact that $  \int_{T}fg'  $ is
convergent, the integral on the left side is convergent and our formula holds.  
$\: \Box$
\begin{lemma}
\label{Le11}
  Let $ f,g $ be continuous on a segment $ T \subset R^{k} $ and for
a $j$ let $ f_{\mid j},g_{\mid j} $ exist and be continuous on $ \Int T  $ with
$g_{\mid j}$ bounded. Then we have
\begin{eqnarray}
\label{co6}
 \int_{T}f_{\mid j}gdt= \int_{\Tkj}[fg]_{T_{j}}d\tkj - \int_{T}fg_{\mid j}dt.
\end{eqnarray}
\end{lemma}
\par Proof.  For simplicity of notations, let us adopt $ T=<0,1>^{k}$. Assume
that $ \mathcal P $ is 
a finite cover of  $  \Fr T $ and $\mathcal P$ consists of segments. The
mapping $ (fg)_{\mid j} $ is integrable on $ T \setminus \bigcup \mathcal P$ and
we may apply the Fubini theorem.  In consequence we have
$ |\int_{T \setminus \bigcup \mathcal P}(fg)_{\mid j}dt - \int_{\Tkj}[fg]_{T_{j}}d\tkj| \leq
\sum\{\int_{\Pkj}|[(fg)(\cdot,\tkj)]_{P_{j}}|d\tkj\!: P \in \mathcal P\}$.
The mapping $ fg $ is uniformly continuous and therefore the sum on the right side of our inequality is small
for $ P_{j} $ containing $0$ or $1$ (such $P$ are ``thin'').  On the other
hand $P$ which do not contain $ (0,\tkj) $ nor $ (1,\tkj) $ are ``flat'' (the measure
of $ \bigcup \Pkj $ in $ \Tkj $ is small).  Therefore the right side of the previous
inequality converges to zero, i.e. the following equality holds
\[
\qquad \int_{T}(fg)_{\mid j}dt=\int_{\Tkj}[fg]_{T_{j}}d\tkj
\]
and the integral on the left side is convergent.
On the other hand $  \int_{T}fg_{\mid j}dt  $ is convergent and we obtain (\ref{co6}). 
$\: \Box$
\par For simplicity let us adopt $ \stackrel{k-1}{\omega}
=b(y)\bigwedge dy_{I}=b_{i_{1},\ldots,i_{k-1}}(y)
dy_{i_{1}} \wedge\cdots\wedge d_{i_{k-1}}$. Clearly it is sufficient to
prove the Stokes formula for such a form.
\par As a direct consequence of Lemma \ref{Le11} we obtain
\begin{lemma}
\label{Le12}
  For $ T=<0,1>^{k} $ let us assume that for $ e\!: T \rightarrow R^{n} $
the partial derivatives of $e$ (one-sided on $  \Fr T$) are continuous and
$  b \circ e \!: T \rightarrow R $ is continuous.  If $  (b \circ e)_{\mid j}$, $ \deteIj  $ 
are continuous on $ \Int T $ and $  \deteIj  $ are bounded, then
\begin{eqnarray}
\label{co7}
\int_{<0,1>^{k}}(b\circ e)_{\mid j}(t)
det \frac{De_{I}}{Dt_{\stackrel{-}{j}}}dt=
\int_{<0,1>^{k-1}} [(b\circ e)(t)
det \frac{De_{I}}{Dt_{\stackrel{-}{j}}}]_{t_{j}=0}^{1} dt_{\stackrel{-}{j}}\\ \nonumber
-\int_{<0,1>^{k}}(b\circ e)(t)(det \frac{De_{I}}{Dt_{\stackrel{-}{j}}})_{\mid j}dt
\end{eqnarray}
holds.
\end{lemma}
\par From the properties of determinants follows
\begin{eqnarray*}
\sum_{i=1}^{n}b_{\mid i}(y)\mid_{y=e(t)}
det \frac{D(e_{i},e_{I})}{D(t_{1},\ldots,t_{k})}=
\sum_{i=1}^{n}det\left( \begin{array}{c}
\begin{array}{ccc}
b_{\mid i}e_{i \mid 1} & \ldots & b_{\mid i}e_{i \mid k}
\end{array}\\
\frac{De_{I}}{Dt}
\end{array} \right)=\\
det\left( \begin{array}{c}
\begin{array}{ccc}
\sum b_{\mid i}e_{i \mid 1} & \ldots & \sum b_{\mid i}e_{i \mid k}
\end{array}\\
\frac{De_{I}}{Dt}
\end{array} \right)=
\sum_{j=1}^{k}(-1)^{j-1}(b \circ e)_{\mid j} det \frac{De_{I}}{Dt_{\stackrel{-}{j}}}.
\end{eqnarray*}
\par From the previous formula and (\ref{co6}) we obtain
\begin{eqnarray}
\label{co8}
\int_{<0,1>^{k}}\sum_{i=1}^{n}b_{\mid i}(y)\mid_{y=e(t)}
det \frac{D(e_{i},e_{I})}{D(t_{1},\ldots,t_{k})}dt=\\ \nonumber
\int_{<0,1>^{k}}\sum_{j=1}^{k}(-1)^{j-1}(b \circ e)_{\mid j}
    det \frac{De_{I}}{Dt_{\stackrel{-}{j}}}dt=\\ \nonumber
\sum_{j=1}^{k}\int_{<0,1>^{k-1}} (-1)^{j-1}[(b\circ e)(t)
det \frac{De_{I}}{Dt_{\stackrel{-}{j}}}]_{t_{j}=0}^{1} dt_{\stackrel{-}{j}}-\\ \nonumber
\sum_{j=1}^{k}\int_{<0,1>^{k}}(b\circ e)(t)(-1)^{j-1}(det \frac{De_{I}}{Dt_{\stackrel{-}{j}}})_{\mid j}dt=\\ \nonumber
\sum_{j=1}^{k}\int_{<0,1>^{k-1}} (-1)^{j-1}[(b\circ e)(t)
det \frac{De_{I}}{Dt_{\stackrel{-}{j}}}]_{t_{j}=0}^{1} dt_{\stackrel{-}{j}}-\\
\nonumber \int_{<0,1>^{k}}(b \circ e)(t) det B(t) dt,
\end{eqnarray}
where
\begin{eqnarray}
\nonumber det B=det\left( \begin{array}{c}
\begin{array}{ccc}
\frac{\partial}{\partial t_{1}} & \ldots & \frac{\partial}{\partial t_{k}}
\end{array}\\
\frac{De_{I}}{Dt}
\end{array} \right).
\end{eqnarray}
The last integral is to be calculated as the sum of integrals of partial
derivatives of determinants.
\begin{corollary}
\label{Cor13}
  Under the assumptions of Lemma \ref{Le12} for $  j=1,\ldots,k$, formula
(\ref{co8}) holds.
\end{corollary}
\begin{remark}
\label{Re14}
  Formula (\ref{co8}) is satisfied if on the both sides we renumber
$ t_{1},\ldots,t_{k} $ in the same way (the same permutation of
columns).  Clearly (\ref{co8}) is independent of $ I=i_{1}\ldots i_{k-1}$.
\end{remark}
\par For $e$ being of class $ C^{2} $ we have $ det B=0 $, as the mixed
second order partial derivatives are equal.
\par Now we are going to present another case in which $ det B=0$.
\par For a set $A$ as in Definition \ref{De2} we define a map
$ c\!: <0,1>^{k} \rightarrow R^{k} $ as follows
\begin{eqnarray}
\label{co9}
\left\{\begin{array}{l}
c_{1}(t_{1})=(1-t_{1})f_{1}+t_{1}g_{1},\\
c_{2}(t_{1},t_{2})=(1-t_{2})f_{2}(c_{1}(t_{1}))+t_{2}g_{2}(c_{1}(t_{1})),\\
\vdots\\
c_{k}(t_{1},\ldots,t_{k})=(1-t_{k})f_{k}(c_{1},\ldots,c_{k-1})+
            t_{k}g_{k}(c_{1},\ldots,c_{k-1}).
\end{array}\right.
\end{eqnarray}
\begin{lemma}
\label{Le15}
  For $ e=c $ (see (\ref{co9}))  we have $ det B=0 $  (see (\ref{co8}))
in points where $c$ is of class $ C^{1}$.
\end{lemma}
\par Proof.  The matrix $ \frac{Dc}{Dt} $ contains only zeros over the main
diagonal, and therefore matrix $B$ has the same property (excluding the
first row). On the diagonal we have $ (g_{j}-f_{j})(c_{1},\ldots,c_{j-1})$ with continuous first order
partial derivatives. The ``determinant operator'' is the sum of the elements of the first row multiplied
by determinants with zeros over the main diagonal (and by $1$ or $-1$). 
Hence the mixed second order partial derivatives are continuous, which means equal
and $ det B=0$. 
$\: \Box$
\par The next step is the following
\begin{lemma}
\label{Le16}
  Let $ A \subset R^{k} $ be a set and let $  \PhA  $ be of class $ C^{2}$.  Then for 
$  e= \Phi \circ c $ (see (\ref{co9})) $ det B=0 $ holds in points where $c$ is of class $ C^{1} $.
\end{lemma}
\par Proof.  The determinant
\begin{eqnarray*}
det\left( \begin{array}{c}
\begin{array}{ccc}
\frac{\partial}{\partial t_{1}} & \ldots & \frac{\partial}{\partial t_{k}}
\end{array}\\
\frac{D(\Phi \circ c)_{I}}{Dt}
\end{array} \right)
\end{eqnarray*}
is the sum of determinants of the form
\begin{eqnarray*}
det\left(
\begin{array}{ccc}
\frac{\partial}{\partial t_{1}} & \ldots & \frac{\partial}{\partial t_{k}}\\
\Phi_{i_{1} \mid i}c_{i \mid 1}  & \ldots & \Phi_{i_{1} \mid i}c_{i \mid k}\\
\Phi_{i_{k-1} \mid j}c_{j \mid 1}  & \ldots & \Phi_{i_{k-1} \mid j}c_{j \mid k}\\
\end{array} \right).
\end{eqnarray*}
While differentiating the rows we are interested in checking the derivatives
of $ c_{i \mid j} $ which can be not continuous.  Therefore
$ \Phi_{i_{m} \mid i} $ can be excluded and finally we have the problem
of continuity of
\begin{eqnarray*}
det\left(
\begin{array}{ccc}
\frac{\partial}{\partial t_{1}} & \ldots & \frac{\partial}{\partial t_{k}}\\
c_{i \mid 1}  & \ldots & c_{i \mid k}\\
c_{j \mid 1}  & \ldots & c_{j \mid k}\\
\end{array} \right).
\end{eqnarray*}
If some rows are identical, this determinant equals to zero, if not - we
apply Lemma \ref{Le15}.  Hence the mixed derivatives are continuous which means
equal and $  det B=0$. 
$\: \Box$
\par Now we put some facts together.
\begin{lemma}
\label{Le17}
  Let $ A \subset R^{k} $ be a normal set and let $  \PhA  $ be of class $ C^{2} $ (or class 
$ C^{1} $ and one to one).
If $ b\!: \Phi (A) \rightarrow R  $ is continuous, 
$  b_{\mid j} $ are continuous and bounded on
$ \Phi(\Int A)  $ then for $c$ given by (\ref{co9})
\begin{eqnarray}
\label{co10}
\int_{<0,1>^{k}} \sum_{i=1}^{n}b_{\mid i}(y)\mid_{y=(\Phi \circ c)(t)} 
det \frac{D((\Phi \circ c)_{i}(t),(\Phi \circ c)_{I}(t))}{D(t_{1},\ldots,t_{k})} dt =\\ \nonumber
\sum_{j=1}^{k}\int_{<0,1>^{k-1}} (-1)^{j-1}[(b\circ \Phi \circ c)(t)
det \frac{D(\Phi \circ c)_{I}(t)}{D\tkj}]_{t_{j}=0}^{1} d\tkj
\end{eqnarray}
holds.
\end{lemma}
\par Proof.  We can disregard $  \Fr A \setminus \Fr \Int A $, 
as there we have zeros on both sides of (\ref{co10}).  The
$  C^{2}  $ case is an immediate consequence
of Corollary \ref{Cor13} and Lemma \ref{Le16}.  Now let $A$ be
normal and $\Phi$ of class $ C^{1} $ and one to one.  In such a case
we consider $ f_{i}^{\epsilon}=min\{f_{i}+\epsilon,g_{i}\}$, $ g_{i}^{\epsilon}=
\mx\{g_{i}-\epsilon,f_{i}^{\epsilon}\} $ and the set $ D_{\epsilon} $ consisting of
points satisfying (\ref{co1}) for $ f_{i}^{\epsilon},g_{i}^{\epsilon} $
in place of $ f_{i},g_{i}$, $ i=1,\ldots,k$.  The set $ A_{\epsilon}=
\overline{\Int D_{\epsilon}} $ is a countable sum of normal sets (the points
with $ f_{i}^{\epsilon}=g_{i}^{\epsilon} $ can be disregarded, as the
integrals are reduced to zero there - see (\ref{co10})). For any open set
$ U \subset R^{k} $ such that $ A_{\epsilon} \subset U \subset A $ there
exist $ \Phi_{n} $ of class $ C^{\infty} $ converging uniformly  with the first derivatives on
$ A_{\epsilon} $ to $\Phi$ and such that $ \Phi_{n}(A_{\epsilon}) \subset \Phi_{n}(U) \subset
\Phi_{n}(\Int A) \subset \Int \Phi(A) $ (in $ \Phi(A)$), as $\Phi$ is one
to one.  Therefore $ b \circ \Phi_{n} $ is well defined.
Now (\ref{co10}) holds for $ \Phi_{n} $ and $ A_{\epsilon} $.  Clearly, from the
uniform convergence with the first order derivatives, we obtain (\ref{co10})
for $\Phi$ on $ A_{\epsilon}$.  We have $  A_{\epsilon}  \rightarrow  A  $ (by measure)
and hence (\ref{co10}) holds for $A$ in place of $ A_{\epsilon}$, as the right side is integrable. 
$\: \Box$
\begin{remark}
\label{Re18}
  If in Lemma \ref{Le17} we assume that for a neighbourhood $V$ of $ \Phi(A)$,
$  b\!: V \rightarrow R  $  is continuous and
$  b_{\mid i} $ are continuous and bounded, then it is sufficient for $  \PhA  $
to be of class $ C^{1}$.  In the respective proof we consider $ \Phi_{n} $ with
$ \Phi_{n}(A_{\epsilon}) \subset V$.
\end{remark}
\par Now Theorem \ref{Th1} is almost proved ;).  All we need to complete
the proof is to get rid of parametrization in (\ref{co10}).
\begin{lemma}
\label{Le19}
  Let $ A \subset R^{k} $ be a normal set and let $  \PhA  $ be of class $ C^{1}$. Assume that
$ b\!: \Phi (A) \rightarrow R  $ is continuous, 
$  b_{\mid j} $ are continuous and  bouned on
$ \Phi(\Int A)$.  Then
\begin{eqnarray}
\label{co11}
\int_{<0,1>^{k}} \sum_{i=1}^{n}b_{\mid i}(y)\mid_{y=(\Phi \circ c)(t)} 
det \frac{D((\Phi \circ c)_{i}(t),(\Phi \circ c)_{I}(t))}{D(t_{1},\ldots,t_{k})} dt =\\ \nonumber
\int_{A} \sum_{i=1}^{n}b_{\mid i}(y)\mid_{y=\Phi(x)} det \frac{D(\Phi_{i}(x),\Phi_{I}(x))}{Dx}dx
\end{eqnarray}
holds.
\end{lemma}
\par Proof. We have $ b_{\mid i}(y)\mid_{y=(\Phi \circ c)(t)} 
det \frac{D((\Phi \circ c)_{i}(t),(\Phi \circ c)_{I}(t))}{D(t_{1},\ldots,t_{k})}=\\
b_{\mid i}(y)\mid_{y=\Phi(x)} det \frac{D(\Phi_{i}(x),\Phi_{I}(x))}{Dx}_{\mid x=
c(t)} \cdot det \frac{Dc(t)}{Dt}$. In view of Remark \ref{Re14}
we may assume
$ det \frac{Dc(t)}{Dt}=(g_{1}-f_{1}) \cdot, \ldots, \cdot(g_{k}-f_{k}) \geq 0$.
This determinant equals to zero on a closed set $B$ contained in
$ c^{-1}(\Fr A)$.  Therefore the integral on the left side of (\ref{co11}) can
be restricted to $ M=<0,1>^{k} \setminus B$.  On the other hand $  \Fr A $ is
of measure zero in $ R^{k} $ and it is meaningless for the integral on the 
right side.  The mapping $ c_{\mid M} $ is one to one and it is of class $  C^{1} $
and bounded outside a compact set of measure zero.  According to Lemma \ref{Le17}
the left side of (\ref{co11}) is 
convergent on $M$ and it is sufficient to consider compact subsets of $ \Int A $
being finite sums of segments.  We apply there the change of variables theorem.
$\: \Box$
\par The right side of (\ref{co10}) is more interesting.  Let us consider the 
following notations $  G_{j}(\xkj)=(g_{j}(\xkj),\xkj)$, $  F_{j}(\xkj)=(f_{j}(\xkj),\xkj) $
and $  \Psi(\xkj)=(\Phi \circ G_{j})(\xkj)$.  Then we have
$ det \frac{D(\Phi \circ c)_{I}(t)}{D\tkj}]_{t_{j}=1}=
det \frac{D(\Psi \circ \ckj)_{I}}{D\tkj}(1,\tkj)=
det \frac {\Psi_{I}(\xkj)}{D\xkj}_{\mid x_{j}=\ckj(1,\tkj)} \cdot 
det \frac{D \ckj}{D \tkj}(1,\tkj)$.  Hence we obtain
\begin{eqnarray*}
H:=\int_{<0,1>^{k-1}} (-1)^{j-1}[(b\circ \Phi \circ c)(t)
det \frac{D(\Phi \circ c)_{I}(t)}{D\tkj}]_{\mid t_{j}=0}^{1} d\tkj=\\
\int_{<0,1>^{k-1}} (-1)^{j-1}(b\circ \Phi \circ c)(1,\tkj)
det \frac{D(\Phi \circ G_{j})_{I}(\xkj)}{D\xkj}_{\mid \xkj=\ckj(1,\tkj)} \cdot
det \frac{D \ckj(1,\tkj)}{D \tkj}d\tkj -\\
\int_{<0,1>^{k-1}} (-1)^{j-1}(b\circ \Phi \circ c)(0,\tkj)
det \frac{D(\Phi \circ F_{j})_{I}(\xkj)}{D\xkj}_{\mid \xkj=\ckj(0,\tkj)} \cdot
det \frac{D \ckj(0,\tkj)}{D \tkj}d\tkj.
\end{eqnarray*}
We have $  \ckj\!: R^{k-1} \rightarrow R^{k-1} $ and like in the proof
of Lemma \ref{Le19} it is sufficient to consider the set of points $ \tkj $
with  $ det \frac{D \ckj}{D \tkj} \neq 0$.  In addition the integrals (on the right side)
on $ (g_{j}-f_{j})^{-1}(0) $ reduce each other and finally we are interested in 
points $  \xkj  $ for which $  G_{j}(\xkj) \neq F_{j}(\xkj) $ takes place.  By the 
change of variables theorem we obtain
\begin{eqnarray}
\label{co12}
H=\int_{S_{j}}[ (b\circ \Phi \circ G_{j})(\xkj) \cdot (-1)^{j-1}
det \frac{D(\Phi \circ G_{j})_{I}(\xkj)}{D\xkj} -\\
\nonumber (b\circ \Phi \circ F_{j})(\xkj) \cdot (-1)^{j-1}
det \frac{D(\Phi \circ F_{j})_{I}(\xkj)}{D\xkj}] d\xkj
\end{eqnarray}
where $  S_{j}=G^{-1}_{j}(\Fr \Int A) $ (condition
$  det \frac{D\ckj}{D\tkj}>0 $ holds there outside a compact set of ($k-1$)-dimensional
measure zero).  The integrals in condition (\ref{co12}) can be transformed
to a more compact form.
We have 
\begin{eqnarray*}
(-1)^{j-1}det \frac{D(\Phi \circ G_{j})_{I}(\xkj)}{D\xkj} =\\
det \left( \begin{array}{ccccc}
0 & \ldots & 1 & \ldots & 0\\
(\Phi \circ G_{j})_{I \mid 1}(\xkj) & \ldots & \Phi_{I \mid j}(z) & \ldots & (\Phi \circ G_{j})_{I \mid k}(\xkj)
\end{array}\right)=\\
det \left( \begin{array}{ccccc}
-g_{j\mid 1(\xkj)} & \ldots & 1 & \ldots & -g_{j\mid k(\xkj)}\\
(\Phi_{I \mid 1}(z) & \ldots & \Phi_{I \mid j}(z) & \ldots & (\Phi_{I \mid k}(z)
\end{array}\right)
\end{eqnarray*}
for $ z=G_{j}(\xkj) $, as $  (\Phi \circ G_{j})_{\mid s}(\xkj)=
\Phi_{\mid s}(z)+\Phi_{\mid j}(z)g_{j\mid s}(\xkj)$.
The last determinant is obtained  from the second one by multiplying the column $j$ by
$  g_{j\mid s} $ and then by substracting it from the column $s$.
In the first row we have a vector $ \stackrel {\rightarrow}{N}$ which is normal
to $  \Fr A  $ and external to $A$.  Now it is clear that for
\begin{eqnarray}
\label{co13}
W_{\Phi_{I}(z)}=det
\left( \begin{array}{c}
\stackrel{\rightarrow}{n}(z)\\
\frac{D\Phi_{I}(z)}{Dz}
\end{array}\right)
\end{eqnarray}
we have $ H=\int_{L_{j}}(G \circ \Phi)(z) W_{\Phi_{I}}(z) dL_{j} $ for 
$  L_{j}=G_{j}(A_{j}) \cup F_{j}(A_{j}) \subset \Fr \Int A$.  By applying the sum on
$j$ we obtain $  \int_{L} (G \circ \Phi)(z) W_{\Phi_{I}}(z) dL $ ($L=\Fr \Int A$)
in place of the right side in (\ref{co10}).
\begin{theorem}
\label{Th20}
Let $ A \subset R^{k} $ $(k>1) $ be a regular set (Def. \ref{De3}) consisting
of normal sets $ A_{1},\ldots,A_{m} $ and let $  \PhA  $
be of class $ C^{2} $ (or of class $ C^{1} $ and one to one).
 If  $ \stackrel{k-1}{\omega}=\sum b_{I} \bigwedge dy_{I} $ ($ k>1 $) is
a continuous form on $ S=\Phi(A)$, $ b_{I\mid i} $ are continuous and 
bounded on $ \Phi(\bigcup\{\Int A_{j}: j=1,\ldots,m\})$,  then for
$ W_{\Phi_{I}(z)} $ as in (\ref{co13}) and $ L=\Fr \Int A $ we have
\begin{eqnarray}
\label{co14}
\int_{A}\sum_{I}\sum_{i=1}^{n}b_{I \mid i}(y)\mid_{y=(\Phi(x)}
det \frac{D(\Phi_{i}(x),\Phi_{I}(x))}{Dx}dx=\\
\int_{L} \sum_{I} (b_{I} \circ \Phi)(z) W_{\Phi_{I}}(z)dL.\nonumber
\end{eqnarray}
\end{theorem}
\par Proof. In view of Lemmas \ref{Le19}, \ref{Le17} and the reasoning preceding our
theorem, the formula holds for each of the sets $ A_{j}$.  Clearly the integral on the left
side is the sum of integrals on $ A_{j}$.  As regards the right side, the integrals on
$  \Fr \Int A_{i} \cap \Fr \Int A_{j} $ reduce, as the versors $\stackrel{\rightarrow}{n}$
for them are opposite. 
$\: \Box$
\par Both sides of (\ref{co14}) depend on $\Phi$.  In order to obtain Theorem \ref{Th1}
we need one step more.
\begin{lemma}
\label{Le21}
  Let us assume that $ A,B \subset R^{k} $ are compact, $  \PhA$, $  \PsA  $ are continuous,
one to one and of class $ C^{1} $ with bounded derivatives outside compact sets of measure zero, 
$ S=\Phi(A)=\Psi(B) $ and let $  det(\Psi^{-1} \circ \Phi)' \geq 0 $
on $A$ outside a set of measure zero.  Then for $ x,s $ such that
$  \Phi(x)=\Psi(s) $ we have
\begin{eqnarray}
\label{co15}
\frac{det \frac{D\Phi_{I}(x)}{Dx}}{|\Phi'(x)|}=\frac{det \frac{D\Psi_{I}(s)}{Ds}}{|\Psi'(s)|}, I=i_{1},\ldots,i_{k}.
\end{eqnarray}
If in addition   $ \stackrel{k-1}{\omega}=\sum b_{I} \bigwedge dy_{I} $
is a continuous $(k-1)$-form on $S$ then
\begin{eqnarray}
\label{co16}
\int_{A}\sum (b_{I} \circ \Phi)(x) det \frac{D\Phi_{I}(x)}{Dx}dx=
\int_{B}\sum (b_{I} \circ \Psi)(s) det \frac{D\Psi_{I}(x)}{Ds}ds
\end{eqnarray}
holds.
\end{lemma}
\par Proof. We have
\begin{eqnarray*}
det \frac{D\Phi_{I}(x)}{Dx}=
det \frac{D\Psi_{I}(s)}{Ds}_{\mid s=(\Psi^{-1}\circ \Phi)(x)} \cdot det (\Psi^{-1}\circ \Phi)'(x)
\end{eqnarray*}
and
\begin{eqnarray*}
|\Phi'(x)|=|(\Psi \circ \Psi^{-1} \circ \Phi)'(x)|=
|\Psi'(s)|_{\mid s=(\Psi^{-1}\circ \Phi)(x)} \cdot det(\Psi^{-1}\circ \Phi)'(x)
\end{eqnarray*}
the last one factor being not less than zero.  Hence we obtain formula (\ref{co15})). 
 The mappings $\Phi$, $\Psi$ are homeomorphisms on compact sets
$A$, $B$ respectively.  In view of the Sard theorem the set
$ \{x \in A\!: |\Phi'(x)|=0\} $ is of measure zero.  The functions
\begin{eqnarray*}
\frac{det \frac{D\Phi_{I}(x)}{Dx}}{|\Phi'(x)|},\frac{det \frac{D\Psi_{I}(s)}{Ds}}{|\Psi'(s)|}
\end{eqnarray*}
are bounded and continuous outside a set of measure zero and now we obtain
the second formula by applying the change of variables theorem. 
$\: \Box$
\par Condition (\ref{co16}) applies also to the right side of (\ref{co12}) for
$ \Phi \circ G_{j}$, $S_{j}$, $k-1 $ in place of $\Phi$, $A$, $k$ respectively.
Now by applying Theorem \ref{Th20} and Lemma \ref{Le21} (see conditions
(\ref{co10}), (\ref{co11}), (\ref{co12}) and (\ref{co16})) we obtain our Theorem \ref{Th1}.
\vspace{.2in}
\par
\flushleft{Dear Reader,}
\par \hspace{.3in} maybe you are disappointed by seeing no references.
I was disappointed by checking well known versions of the Stokes theorem and finding 
none similar to the one presented here. I was trying to publish this paper for several years. 
Let me know if Theorem \ref{Th1} can be found somewhere. 
\par \hspace{2.7in} Lech Pasicki
\vspace{.5in}
\par
\mbox{Faculty of Applied Mathematics} \linebreak
\mbox{AGH University of Science and Technology} \linebreak
\mbox{Al. Mickiewicza 30} \linebreak
\mbox{30-059 KRAK\'OW, POLAND} \linebreak
\mbox{E-mail: pasicki@agh.edu.pl}

\end{document}